\documentclass{article}[14pt]
\usepackage{amssymb,amsmath,graphicx,enumerate}
\usepackage{longtable,mathrsfs}
\usepackage{hyperref}
\graphicspath {G:\Pachpatte}

\DeclareMathAlphabet{\mathpzc}{OT1}{pzc}{m}{it}

\numberwithin{equation}{section}

\newtheorem{theorem}{Theorem}[section]

\newtheorem{cor}{Corollary}[section]

\newtheorem{definition}{Definition}[section]
\newtheorem{example}{Example}[section]
\newtheorem{remark}{Remark}[section]
\newtheorem{prop}{Proposition}[section]

\begin{document}
\begin{large}
\begin{center}
\large {Fuzzy Hermite-Hadamard type inequality for $r$-preinvex and $(\alpha,m)$-preinvex function}
\end{center}
\begin{center}
                 $^{1}$  Kavita U. Shinde  and  $^{2}$  Deepak B. Pachpatte
\end{center}

\begin{center}
$^{1}$ Department of Mathematics,\\

Dr.Babasaheb Ambedkar Marathwada

University, Aurangabad-431 004 (M.S) India.\\

E-mail:kansurkar14@gmail.com\\

$^{2}$ Department of Mathematics,\\

Dr.Babasaheb Ambedkar Marathwada

University, Aurangabad-431 004 (M.S) India.\\

E-mail: pachpatte@gmail.com\\

\end{center}

\begin{abstract}
  The purpose of this paper is to study the Hermite-Hadamard type inequality for $r$-preinvex and $(\alpha,m)$-preinvex function  which is based  on  Sugeno integral.
\end{abstract}


\section{Introduction}

In 1974 Sugeno stared study of theory of fuzzy integral \cite{sugeno}. It is very useful tool for various application in theoretical and applied statistics which is based on non-additive measures.
  Hanson \cite{han} introduced a generalization of convex functions in terms of invex function.  In   \cite{mnoor}, \cite{vari}, \cite{tweir} authors studied the basic properties and role of preinvex functions in optimization, equilibrium problems and variational inequalities.

 The study of inequalities for Sugeno integral was initiated by Roman-Flores et.al.,  \cite{florom, rom}. H. Agahi et. al., \cite{agahih}, \cite{haga} proved the general Barnes-Godunova-Levin and new general extensions of Chebyshev type inequalities for Sugeno integrals. J. Caballero and K. Sadarangani \cite{g},\cite{h}  studied the Cauchy-Schwarz type inequality  and Chebyshev type inequality for Sugeno integrals.   N. Okur \cite{okur} proved the Hermite-Hadamard type inequality for log-preinvex function using Sugeno integral.   In \cite{HH,kavi,kshinde,abba}
  authors proved the fuzzy Hermite-Hadamard type inequality for convex functions.

 Motivated by  above results in this paper we obtain the Hermite-Hadamard type inequality for $r$-preinvex function and
  $(\alpha,m)$-preinvex function with respect to Sugeno integral.

    \section{Preliminaries}
\qquad\, In this section we give some basic definitions and properties of  Sugeno integral,
 \cite{sugeno},  \cite{1p}, \cite{1z}.

Suppose that $\wp$ is a $\sigma$-algebra of subsets of $X$ and $\mu:\wp\longrightarrow[0,\infty)$ be a non-negative, extended real valued set function.
We say that $\mu$ is a fuzzy measure if it satisfies:
\begin{enumerate}
  \item $\mu(\phi)=0$.
  \item $ E,F\in\wp$ and $E\subset F$ imply $\mu(E)\leq\mu(F)$.
  \item  $\{E_{n}\}\subset\wp, E_{1}\subset E_{2}\subset..., $ imply $\lim_{n\longrightarrow\infty}\mu(E_{n})=\mu(\bigcup_{n=1}^{\infty}E_{n})$.
  \item $\{E_{n}\}\subset\wp$, $E_{1}\supset E_{2}\supset...,$ $\mu(E_{1})<\infty$, imply $\lim_{n\longrightarrow\infty}\mu(E_{n})=\mu(\bigcap_{n=1}^{\infty}E_{n})$.
\end{enumerate}

\par If $f$ is  non-negative real-valued function defined on $X$, we denote the set $\{x\in X:f(x)\geq\alpha\}=\{x\in X:f\geq\alpha\}$
 by $F_{\alpha}$  for $\alpha\geq 0$. Note that if $\alpha\leq\beta$ then $F_{\beta}\subset F_{\alpha}$.

 Let $(X,\wp,\mu)$ be a fuzzy measure space, we denote $M^{+}$ the set of all non-negative measurable functions with respect to $\wp$.
 \begin{definition}(Sugeno \cite{sugeno}).
   Let $(X,\wp,\mu)$ be a fuzzy measure space, $f\in M^{+}$ and $A\in\wp$, the Sugeno integral (or fuzzy integral)
of $f$ on $A$, with respect to the fuzzy measure $\mu$, is defined as

 $$ (s)\int_{A}f d\mu=\underset{\alpha\geq0}{\bigvee}[\alpha\wedge\mu(A\cap F_{\alpha})],$$
when $A=X$, \
$$ (s)\int_{X}f d\mu=\underset{\alpha\geq0}{\bigvee}[\alpha\wedge\mu(F_{\alpha})],$$
where $\bigvee$ and $\wedge$ denote the operations sup and inf on $[0,\infty)$, respectively.
 \end{definition}
The properties of Sugeno integral are  well known and can be found in \cite{2z} as follows.
\begin{prop}\label{z2}
   Let $(X,\wp,\mu)$ be fuzzy measure space, $A,B\in \wp$ and $f,g\in M^{+}$ then:
\begin{enumerate}
  \item  $(s)\int_{A} f d \mu \leq \mu(A)$.
  \item  $(s)\int_{A} kd \mu = k\wedge\mu(A)$, $k$ for non-negative constant.
  \item $(s)\int_{A} f d \mu \leq (s)\int_{A}gd\mu$, for $f\leq g$.
  \item $\mu(A\cap\{f\geq\alpha\})\geq\alpha\Longrightarrow(s)\int_{A}fd\mu\geq\alpha$.
  \item$\mu(A\cap\{f\geq\alpha\})\leq\alpha\Longrightarrow(s)\int_{A}fd\mu\leq\alpha$.
  \item$(s)\int_{A} f d \mu >\alpha \Longleftrightarrow$ there exists $\gamma>\alpha$ such that $\mu(A\cap\{f\geq\gamma\})>\alpha$.
  \item $(s)\int_{A} f d \mu <\alpha \Longleftrightarrow$ there exists $\gamma<\alpha$ such that $\mu(A\cap\{f\geq\gamma\})<\alpha$.
\end{enumerate}
\end{prop}
\begin{remark}\label{z1}
  Consider the distribution function $F$ associated to $f$ on $A$, that is, $F(\alpha)=\mu(A\cap\{f\geq\alpha\})$.Then due to
$(4)$ and $(5)$ of Proposition \ref{z2},  we have $F(\alpha)=\alpha\Longrightarrow (s)\int_{A}fd\mu=\alpha$. Thus, from a numerical point of view, the fuzzy integral can be calculated solving the equation $F(\alpha)=\alpha$.
\end{remark}

 Let $f:K\longrightarrow\mathbb{R}$  and
$\eta(.,.):K \times K\longrightarrow \mathbb{R}^{n}$
be continuous functions, where $K\subset\mathbb{R}^{n}$ is a nonempty closed set.

\quad We use the notations, $\langle.,.\rangle$ and $\|.\|$, for
inner product and norm, respectively.    Now we give some definitions  and condition which will be useful  subsequent discussion.

\begin{definition}\label{get}(\cite{mnoor,tweir})
	Let $u\in K $. Then, the set $K$ is said to be invex at $u\in K$ with respect to $\eta(.,.)$ if
	\begin{equation}\label{we1}
	u+t\eta(v,u)\in K, \quad \forall u,v\in K, \, t\in[0,1].
	\end{equation}
	The invex set $K$ is also called a $\eta$
	connected set.
\end{definition}

\quad Invex set has a clear geometric interpretation which
says that there is a path starting from a point $u$ which is contained in $K$.
It does  not require that the point $v$ should be one of the end points of the path. This observation
plays an important role \cite{tant}.

\quad If   $v$ is an end point
of the path for every pair of points, $u,v \in K,$ then $\eta(v,u)=v-u$ and consequently invexity reduces
to convexity. Thus every convex set is also an invex set with respect to  $\eta(v,u)=v-u,$ but its
converse  is not necessarily true \cite{tweir}, \cite{myang}.

\begin{definition}(T. Weir \cite{tweir}).
	The function $f$ on the invex set $K$ is said to be preinvex with respect to $\eta$ if
	\begin{equation}\label{we2}
	f(u+t\eta(v,u))\leq (1-t)f(u)+t f(v), \quad \forall u,v\in K, \, t\in[0,1].
	\end{equation}
\end{definition}

In \cite{srmohan}, Mohan and Neogy has given the condition for function $\eta$ known as, \\
\textbf{Condition $C$} :  Let $K\subseteq \mathbb{R}$ be an open invex subset to $\eta: K\times K\longrightarrow\mathbb{R}.$
For any $x,y\in K$ and any $t\in[0,1],$
\begin{align*}
\eta(y, y+t\eta(x,y))=&-t\eta(x,y),\\
\eta(x,y+t\eta(x,y))=& (1-t) \eta(x,y).
\end{align*}
If for any $x,y\in K$ and $t_1, t_2\in [ 0,1],$ we have
\begin{equation*}
\eta(y+t_2\eta(x,y),y+t_1\eta(x,y))=(t_2-t_1) \eta(x,y).
\end{equation*}

\quad The concepts of the invex and preinvex functions have played very important roles in the development of generalized convex programming,
see  \cite{Pachpatte, MAnoor,NoorMA,srmohan,NM,pini}.

\quad In \cite{antcz,antczak},  Antczak introduced  the concept of $r$-invex and $r$-preinvex functions  which as follows.
\begin{definition}
	A positive function $f$ on the invex set $K$ is said to  be $r$-preinvex  with respect to $\eta$ if, for each
	$u,v\in K,$     $t\in [0,1]$.
	\begin{equation*}
	f(u+t\eta(v,u))\leq \begin{cases}
	((1-t)f^r(u)+tf^r(v))^{1/r} , & r\neq 0   \\
	(f(u))^{1-t}(f(v))^t, &  r=0.
	\end{cases}
	\end{equation*}
	Note that $0$-preinvex functions are logarithmic  preinvex  and $1$-preinvex functions are classical preinvex functions.
	If $f$ is $r$-preinvex function then $f^r$ is preinvex function $(r>0).$
\end{definition}
\quad The $m$-preinvex function is defined as
\begin{definition} \label{d1}
	\cite{latif}
	The function $f$ on the invex set $K\subseteq [0,b^{*}]$, $b^{*}>0$ is said to be $m$-preinvex with respect to $\eta$ if
	\begin{equation*}
	f(u+t\eta(v,u))\leq (1-t)f(u)+mtf(\frac{v}{m}),
	\end{equation*}
	holds for all $u,v\in K$, $t\in [0,1]$ and $m\in(0,1]$. The function $f$ is said to be $m$-preconcave if and only if $-f$ is $m$-preinvex.
\end{definition}

\begin{definition} \label{d2} \cite{latif}
	The function $f$ on the invex set $K\subseteq [0,b^{*}]$,  $b^{*}>0$ is said to be $(\alpha,m)$-preinvex function with respect to $\eta$ if
	\begin{equation*}
	f(u+t\eta(v,u))\leq (1-t^{\alpha})f(u)+mt^{\alpha}f(\frac{v}{m}),
	\end{equation*}
	holds for all $u,v \in K$, $t\in[0,1]$ and $(\alpha,m)\in(0,1]\times(0,1]$. The function $f$ is said to be $(\alpha,m)$-preincave if and only if $-f$ is $(\alpha,m)$-preinvex.
\end{definition}
\begin{remark}\cite{latif}
	If we put   $m=1$ in Definition \ref{d1}, then  we get the  definition of preinvexity. If we put $\alpha=m=1$, then Definition \ref{d2}  becomes the definition of preinvex function.  Every $m$-preinvex  function and $(\alpha,m)$-preinvex functions are $m$-convex and $(\alpha,m)$-convex with respect to $\eta(v,u)=v-u $ respectively.
\end{remark}

\section{Fuzzy Hadamard type inequality for $r$-preinvex function}
Now in this section we give results obtained on Hadamard type inequality for $r$-preinvex function with respect to Sugeno integral.

\quad In \cite{wasim} W. Ul-Haq and J. Iqbal proved the following Hermite-Hadamard type inequality for $r$-preinvex function.
\begin{theorem}
	Let $f:K=[a,a+\eta(b,a)]\longrightarrow (0,\infty)$ be an $r$-preinvex function on the interval of real numbers $K^o$ (interior of K) and
	$a,b\in K^o $ with $a< a+ \eta(b,a).$ Then the following inequalities holds
	\begin{equation}\label{pap1}
	\frac{1}{\eta(b,a)}\int_{a}^{a+\eta(b,a)} f(x)dx \leq \bigg[\frac{f^r(a)+f^r(b)}{2}\bigg]^{1/r}, \, \, r\geq 1.
	\end{equation}
\end{theorem}
Now we consider an example.
\begin{example}
	Consider $X=[0,\eta(1,0)]$ and let $\mu$ be the Lebesgue measure on $X$. If we take  $f(x)= \frac{x^4}{2}$  be a  non-negative and $\frac{1}{2}$-preinvex function
	on $[0,\eta(1,0)].$
	From Remark \ref{z1}, we have
	\begin{align}\label{dc1}
	\nonumber  F(\beta) =& \mu\big([0,\eta(1,0)]\cap\big\{x|\frac{x^4}{2}\geq \beta\big\}\big)\\
	\nonumber  =& \mu\big((2\beta)^{1/4},\eta(1,0)\big)\\
	=& \eta(1,0)-(2\beta)^{1/4},
	\end{align}
	and the solution of \eqref{dc1} is
	\begin{equation*}
	\eta(1,0)-(2\beta)^{1/4}=\beta,
	\end{equation*}
	where  $0\leq \eta(1,0)\leq 1$ to $ 0\leq \beta \leq 0.2023,$    we get
	\begin{equation*}
	0\leq (s)\int_{0}^{\eta(1,0)}f d\mu= (s)\int_{0}^{\eta(1,0)}\frac{x^4}{2} d\mu \leq 0.2023.
	\end{equation*}
	\begin{equation}\label{pap2}
	0.2023\leq \frac{1}{\eta(1,0)}(s)\int_{0}^{\eta(1,0)}\frac{x^4}{2} d\mu < \infty.
	\end{equation}
	On the other hand
	\begin{equation}\label{pap3}
	\bigg[\frac{f^{1/2}(0)+f^{1/2}(\eta(1,0))}{2}\bigg]^2= \frac{\eta(1,0)^4}{8}=\frac{1}{8}=0.125,
	\end{equation}
	where $0\leq \eta(1,0)\leq 1.$
	
	\quad From inequalities \eqref{pap2}, \eqref{pap3} and  Boolean operator on real numbers, it is seen that
	Hermite-Hadamard type inequality \eqref{pap1} is not valid in fuzzy context.
\end{example}

\quad Now we give the Hermite-Hadamard type inequality for Sugeno integral with respect to $r$-preinvex function.
\begin{theorem}\label{of1}
	Let $r>0$ and $\mu$ be the Lebesgue  measure on $\mathbb{R}$,  $f:[0,\eta(1,0)]\longrightarrow[0,\infty)$ be $r$-preinvex function with
	$f(0)\neq f(\eta(1,0)).$\\
	\textbf{Case 1.} If $f(\eta(1,0))>f(0),$ then
	\begin{equation*}
	(s)\int_{0}^{\eta(1,0)}f d\mu \leq min\{\beta,\eta(1,0)\},
	\end{equation*}
	where $\beta$ satisfies the following equation
	\begin{equation}\label{pap4}
	(f^r(\eta(1,0))-f^r(0))\beta + \eta(1,0)\beta^r-\eta(1,0)f^r(\eta(1,0))=0.
	\end{equation}
	\textbf{Case 2.} If $f(0)>f(\eta(1,0)),$ then
	\begin{equation*}
	(s)\int_{0}^{\eta(1,0)}f d\mu \leq min\{\beta,\eta(1,0)\},
	\end{equation*}
	where $\beta$ satisfies the following equation
	\begin{equation}\label{pap6}
	(f^r(\eta(1,0))-f^r(0))\beta -\beta^r \eta(1,0)+\eta(1,0) f^r(0)=0.
	\end{equation}
\end{theorem}

\textbf{Proof.}
As $f$ is a $r$-preinvex  function  for $x\in [0,\eta(1,0)]$ we have
\begin{equation*}
f(0+x\eta(\eta(1,0)))\leq ((1-x)f^r(0)+ x f^r(\eta(1,0)))^{1/r},
\end{equation*}
and from Condition C, we have
\begin{equation*}
\eta(\eta(1,0),0)= \eta(0+1.\eta(1,0),0+0.\eta(1,0))=\eta(1,0).
\end{equation*}
Therefore,
\begin{align*}
f(x)=&f\bigg(\frac{x}{\eta(1,0)}.\eta(1,0)\bigg)\\
\leq & \bigg[\bigg(1-\frac{x}{\eta(1,0)}\bigg)f^r(0)+\bigg(\frac{x}{\eta(1,0)}\bigg)f^r(\eta(1,0))\bigg]^{1/r}\\
=& g(x).
\end{align*}
By Proposition \eqref{z2}, we have
\begin{align}\label{pap7}
\nonumber (s)\int_{0}^{\eta(1,0)} f(x) d\mu =& (s)\int_{0}^{\eta(1,0)} f\bigg(\frac{x}{\eta(1,0)}.\eta(1,0)\bigg) d\mu\\
\nonumber \leq & (s)\int_{0}^{\eta(1,0)} \bigg[\bigg(1-\frac{x}{\eta(1,0)}\bigg)f^r(0)\\
\nonumber &\qquad\qquad+\bigg(\frac{x}{\eta(1,0)}\bigg)f^r(\eta(1,0))\bigg]^{1/r} d\mu\\
=& (s)\int_{0}^{\eta(1,0)} g(x) d\mu.
\end{align}
To calculate the right hand side of \eqref{pap7}, we consider the distribution function $F$ given  by
\begin{align}\label{pap8}
\nonumber  F(\beta) =& \mu\bigg([0,\eta(1,0)]\cap \bigg\{x|g(x)\geq \beta\bigg\}\bigg)\\
\nonumber   =& \mu\bigg([0,\eta(1,0)]\cap \bigg\{x|\bigg[\bigg(1-\frac{x}{\eta(1,0)}\bigg)f^r(0)\\
&\qquad\qquad \qquad +\bigg(\frac{x}{\eta(1,0)}\bigg)f^r(\eta(1,0))\bigg]^{1/r}\geq \beta\bigg\}\bigg).
\end{align}
\textbf{Case 1.} If $f(\eta(1,0))>f(0),$ then  from \eqref{pap8}, we have
\begin{align}\label{pap9}
\nonumber  F(\beta)=& \mu\bigg([0,\eta(1,0)]\cap \bigg\{x|x\geq \eta(1,0)\frac{\beta^r-f^r(0)}{f^r(\eta(1,0))-f^r(0)}\bigg\}\bigg)\\
\nonumber   =& \mu\bigg(\eta(1,0)\frac{\beta^r-f^r(0)}{f^r(\eta(1,0))-f^r(0)}, \eta(1,0)\bigg)\\
=& \eta(1,0)- \eta(1,0)\frac{\beta^r-f^r(0)}{f^r(\eta(1,0))-f^r(0)},
\end{align}
and the solution of \eqref{pap9} is $F(\beta)=\beta$,
where $\beta$ satisfies the following equation
\begin{align*}
&\beta(f^r(\eta(1,0))-f^r(0))+\beta^r\eta(1,0) -\eta(1,0) f^r(\eta(1,0))=0.
\end{align*}
By Proposition \ref{z2}  and Remark \ref{z1},  we have
\begin{equation*}
(s)\int_{0}^{\eta(1,0)}f(x) d\mu\leq min\{\beta,\eta(1,0)\}.
\end{equation*}
\textbf{Case 2.} If $f(0)> f(\eta(1,0))$, then from \eqref{pap8}, we have
\begin{align}\label{pap11}
\nonumber F(\beta)=& \mu\bigg([0,\eta(1,0)]\cap \bigg\{x|x\leq\eta(1,0)\frac{\beta^r-f^r(0)}{f^r(\eta(1,0))-f^r(0)}\bigg\}\bigg)\\
\nonumber  =&\mu\bigg(0, \eta(1,0)\frac{\beta^r-f^r(0)}{f^r(\eta(1,0))-f^r(0)}\bigg)\\
=& \eta(1,0)\frac{\beta^r-f^r(0)}{f^r(\eta(1,0))-f^r(0)},
\end{align}
and the solution \eqref{pap11} is $F(\beta)=\beta$,
where $\beta$ satisfies the following equation
\begin{align*}
& \beta (f^r(\eta(1,0))-f^r(0))-\eta(1,0)\beta^r +\eta(1,0)f^r(0)=0.
\end{align*}
By Proposition \ref{z2}  and Remark \ref{z1},  we have
\begin{equation*}
(s)\int_{0}^{\eta(1,0)}f(x) d\mu\leq min\{\beta,\eta(1,0)\}.
\end{equation*}
\begin{example}
	Consider $X=[0,\eta(1,0)]$ and let $\mu$ be the Lebesgue measure on $X$. If we take $f(x)=\frac{x^3}{3}$ be the $ \frac{1}{2}$-preinvex function, where $0\leq \eta(1,0)\leq 1$ from Remark \ref{z2}, we have
	\begin{equation*}
	(s)\int_{0}^{\eta(1,0)} \frac{x^3}{3} d\mu=0.1847.
	\end{equation*}
	From Theorem \ref{of1}, we have
	\begin{equation*}
	0.1847=	(s)\int_{0}^{\eta(1,0)} \frac{x^3}{3} d\mu\leq min\{0.2087,\eta(1,0)\}=0.2087.
	\end{equation*}
\end{example}
\begin{remark}
	In case if we take  $f(0)=f(\eta(1,0))$  in Theorem \eqref{of1}, then we get
	\begin{align*}
	&(s)\int_{0}^{\eta(1,0)} f(x) d\mu \leq (s)\int_{0}^{\eta(1,0)} g(x) d\mu= (s)\int_{0}^{\eta(1,0)} f(0) d\mu\\
	&\qquad\qquad\qquad \qquad\qquad\qquad\qquad\qquad= f(0)\wedge \eta(1,0).
	\end{align*}
\end{remark}
Now we obtain the results for inequality on $r$-preinvex function.
\begin{theorem}\label{of22}
	Let $r> 0$ and $\mu$ be the Lebesgue measure on $\mathbb{R}$,  $f:[a,a+\eta(b,a)]\longrightarrow[0,\infty)$ be
	$r$-preinvex function with $f(a)\neq f(a+\eta(b,a)).$ \\
	\textbf{Case 1.} If $f(a+\eta(b,a))> f(a)$, then
	\begin{equation*}
	(s)\int_{a}^{a+\eta(b,a)} f(x)d\mu \leq min\{\beta,\eta(b,a)\},
	\end{equation*}
	where $\beta$ satisfies the following equation
	\begin{equation}\label{pap13}
	\beta(f^r(a+\eta(b,a))-f^r(a))+\beta^r \eta(b,a)-\eta(b,a)f^r(a+\eta(b,a))=0.
	\end{equation}
	\textbf{Case 2.} If $f(a+\eta(b,a))< f(a), $ then
	\begin{equation*}
	(s)\int_{a}^{a+\eta(b,a)} f(x) d\mu \leq min\{\beta,\eta(b,a)\},
	\end{equation*}
	where $\beta$  satisfies the following equation
	\begin{equation}\label{pap14}
	\beta(f^r(a+\eta(b,a))-f^r(a))-\beta^r\eta(b,a)+\eta(b,a) f^r(a)=0.
	\end{equation}
\end{theorem}

\begin{remark}
	If we take  $f(a)=f(a+\eta(b,a))$ in  Theorem \ref{of22}  we have $g(x)=f(a)$ and by Proposition \ref{z2}, we have
	\begin{align*}
	(s)\int_{a}^{a+\eta(b,a)} f(x) d\mu \leq& (s)\int_{a}^{a+\eta(b,a)} g(x) d\mu\\
	=& (s)\int_{a}^{a+\eta(b,a)} f(a) d\mu\\
	=& f(a)\wedge\mu([a, a+\eta(b,a)]).
	\end{align*}
\end{remark}
\begin{cor}
	Let $r<0,$ and $\mu$ be the Lebesgue  measure on $\mathbb{R}.$ Let $f:[a, a+\eta(b,a)]\longrightarrow [0,\infty)$
	be the $r$-preinvex  with  $f(a)\neq f(a+\eta(b,a)),$ then \\
	\textbf{Case 1.}   If $f(a+\eta(b,a))> f(a),$  then
	\begin{equation*}
	(s)\int_{a}^{a+\eta(b,a)} f(x) d\mu\leq min\{\beta, \eta(b,a)\},
	\end{equation*}
	where $\beta$  satisfies the following equation
	\begin{equation}\label{pap26}
	\beta(f^r(a+\eta(b,a))-f^r(a))-\beta^r \eta(b,a)+\eta(b,a)f^r(a)=0.
	\end{equation}
	\textbf{Case 2.} If $f(a+\eta(b,a)) < f(a),$ then
	\begin{equation*}
	(s)\int_{a}^{a+\eta(b,a)} f(x) d\mu\leq min\{\beta, \eta(b,a)\},
	\end{equation*}
	where $\beta$  satisfies the following equation
	\begin{equation}\label{pap27}
	\beta(f^r(a+\eta(b,a))-f^r(a))+\beta^r\eta(b,a)-\eta(b,a) f^r(a+\eta(b,a))=0.
	\end{equation}
\end{cor}

\section{Fuzzy Hermite-Hadamard type inequality for  $(\alpha,m)$-preinvex function }
Now in this section we give the Hermite-Hadamard type inequality for  $(\alpha,m)$-preinvex function.
In \cite{MAN2}  Noor  proved following the Hermite-Hadamard type inequality for preinvex functions.
\begin{theorem}
	Let $f:[a,a+\eta(b,a)]\rightarrow(0,\infty)$ be a preinvex function on the interval of real numbers $K^\circ$ (the interior of K) and $a,b\in K^\circ$ with $a<a+\eta(b,a)$. Then the following inequality holds,
	\begin{align}\label{or1}
	f\bigg(\frac{2a+\eta(b,a)}{2}\bigg)\leq  \frac{1}{\eta(b,a)}\int_{a}^{a+\eta(b,a)}f(x)dx
	\leq  \frac{f(a)+f(b)}{2}.
	\end{align}
\end{theorem}
Now consider an examples.
\begin{example}
	Consider  $X=[0,\eta(1,0)]$ and let $\mu$ be the Lebesgue measure on $X$. If we take the function $ f(x)= \frac{x^2}{2}$ then
	$f(x)$ is $(1/2,1/3)$-preinvex function.
	From Remark \ref{z2} we have
	\begin{align}\label{or2}
	\nonumber   F(\beta)=& \mu\big([0,\eta(1,0)]\cap\big\{x|\frac{x^2}{2}\geq \beta\big\}\big)\\
	= &\eta(1,0)-\sqrt{2\beta},
	\end{align}
	and the solution of \eqref{or2} is
	\begin{equation*}
	\eta(1,0)-\sqrt{2\beta}=\beta,
	\end{equation*}
	where $0\leq \eta(1,0)\leq 1$ to $0\leq\beta\leq 0.2679$,
	we get
	\begin{equation*}
	0\leq (s)\int_{0}^{\eta(1,0)}fd\mu=(s)\int_{0}^{\eta(1,0)}\frac{x^2}{2} d\mu\leq 0.2679.
	\end{equation*}
	\begin{equation*}
	0.2679 \leq \frac{1}{\eta(1,0)}(s)\int_{0}^{\eta(1,0)}\frac{x^2}{2} d\mu <\infty.
	\end{equation*}
	From right hand side of \eqref{or1} and  for $0\leq \eta(1,0)\leq 1$, we have
	\begin{equation*}
	\frac{f(0)+f(\eta(1,0))}{2}=0.25.
	\end{equation*}
	\quad This proves that the right hand side of \eqref{or1} is not satisfied for Sugeno integral.
\end{example}
\begin{example}
	Consider $X=[0,\eta(1,0)]$ and let $\mu$ be the Lebesgue measure on $X$. If we take the function $f(x)=3x^2$ then
	$f(x)$ is $(1/2,1/3)$-preinvex function. From Remark \ref{z2}, we have
	\begin{align}\label{or3}
	\nonumber   F(\beta)=& \mu\big([0,\eta(1,0)]\cap\big\{x|3x^2\geq \beta\big\}\big)\\
	= &\eta(1,0)-\sqrt{\frac{\beta}{3}},
	\end{align}
	and the solution of \eqref{or3} is
	\begin{equation*}
	\eta(1,0)-\sqrt{\frac{\beta}{3}}=\beta,
	\end{equation*}
	where $0\leq \eta(1,0)\leq 1$ and  $0\leq\beta\leq 0.5657$, we get
	\begin{equation*}
	0\leq (s)\int_{0}^{\eta(1,0)}fd\mu=(s)\int_{0}^{\eta(1,0)}3x^2 d\mu\leq 0.5657.
	\end{equation*}
	\begin{equation*}
	0.5657 \leq \frac{1}{\eta(1,0)}(s)\int_{0}^{\eta(1,0)}3x^2 d\mu <\infty.
	\end{equation*}
	From left hand side of \eqref{or1} and  for $0\leq \eta(1,0)\leq 1$, we have
	\begin{equation*}
	f\bigg(\frac{2a+\eta(b,a)}{2}\bigg)=0.75.
	\end{equation*}
	This proves that the left hand side of \eqref{or1} is not satisfied for Sugeno integral.
\end{example}
 Now in next theorem we prove  Hermite-Hadamard type inequality for Sugeno integral with respect to $(\alpha,m)$-preinvex function.
\begin{theorem}\label{ki5}
	Let $f:[0,\eta(1,0)]\rightarrow[0,\infty)$ be $(\alpha,m)$-preinvex function, $(\alpha,m)\in(0,1)^2$, $f(0)\leq f(\eta(1,0))$ and $\mu$
	be the Lebesgue measure on $\mathbb{R}$. Then
	\begin{equation*}
	(s)\int_{0}^{\eta(1,0)} f(x) d\mu \leq min\{\beta,\eta(1,0)\},
	\end{equation*}
	where $\beta$ satisfies the following equation
	\begin{equation}\label{or6}
	(\eta(1,0)-\beta)^\alpha m f(\frac{\eta(1,0)}{m})-(\eta(1,0)-\beta)^\alpha f(0)- \eta(1,0)^\alpha (\beta-f(0))=0.
	\end{equation}
\end{theorem}
Now we give the Hermite- Hadamard type inequality for decreasing $(\alpha,m)$-preinvex function.
\begin{theorem}\label{ki1}
	Let  $f:[0,\eta(1,0)]\rightarrow[0,\infty)$ be $(\alpha,m)$-preinvex function, $(\alpha,m)\in(0,1]^2$, $f(0)>f(\eta(1,0))$ and $\mu$
	be the Lebesgue measure on $\mathbb{R}$. Then\\
	\textbf{Case 1.} If $m\in (0,\frac{f(\eta(1,0))}{f(0)})$, then
	\begin{equation*}
	(s)\int_{0}^{\eta(1,0)} f(x) d\mu \leq min\{\beta,\eta(1,0)\},
	\end{equation*}
	where $\beta$  satisfies the following equation
	\begin{equation}\label{or14}
	(\eta(1,0)-\beta)^\alpha m f(\frac{\eta(1,0)}{m})-(\eta(1,0)-\beta)^\alpha f(0)- \eta(1,0)^\alpha (\beta-f(0))=0.
	\end{equation}
	\textbf{Case 2.} If $m=\frac{f(\eta(1,0))}{f(0)}$, then
	\begin{equation*}
	(s)\int_{0}^{\eta(1,0)} f(x) d\mu \leq min\{\beta,\eta(1,0)\},
	\end{equation*}
	where $\beta$ satisfies the following equation
	\begin{align}\label{or15}
	\nonumber  & (\eta(1,0)-\beta)^\alpha \frac{f(\eta(1,0))}{f(0)}f\bigg(\frac{f(0)\eta(1,0)}{f(\eta(1,0))}\bigg)\\
	&\quad-(\eta(1,0)-\beta)^\alpha f(0)-\eta(1,0)^\alpha (\beta-f(0))=0.
	\end{align}
	\textbf{Case 3.} If $m\in(\frac{f(\eta(1,0))}{f(0)},\eta(1,0))$, then
	\begin{equation*}
	(s)\int_{0}^{\eta(1,0)} f(x) d\mu \leq min\{\beta,\eta(1,0)\},
	\end{equation*}
	where $\beta$  satisfies the following equation
	\begin{equation}\label{or17}
	\beta^\alpha(m f(\frac{\eta(1,0)}{m})-f(0) )-\eta(1,0)^\alpha (\beta-f(0))=0.
	\end{equation}
\end{theorem}

Now we give the general case of  Theorem \ref{ki5} and \ref{ki1}.
\begin{theorem}\label{ks1}
	Let $f:[a,a+\eta(b,a)]\rightarrow[0,\infty)$ be $(\alpha,m)$-preinvex function, $(\alpha,m)\in(0,1)^2$, $f(a)\leq f(a+\eta(b,a))$ and $\mu$
	be the Lebesgue measure on $\mathbb{R}$. Then
	\begin{equation*}
	(s)\int_{a}^{a+\eta(b,a)} f(x) d\mu \leq min\{\beta,\eta(b,a)\},
	\end{equation*}
	where $\beta$ satisfies the following equation
	\begin{equation}\label{or26}
	(\eta(b,a)-\beta)^\alpha mf(\frac{a+\eta(b,a)}{m})-(\eta(b,a)-\beta)^\alpha f(a)-\eta(b,a)^\alpha (\beta-f(a))=0.
	\end{equation}
\end{theorem}

\begin{theorem}
	Let  $f:[a,a+\eta(b,a)]\rightarrow[0,\infty)$ be $(\alpha,m)$-preinvex function, $(\alpha,m)\in(0,1)^2$, $f(a)> f(a+\eta(b,a))$ and $\mu$
	be the Lebesgue measure on $\mathbb{R}$. Then\\
	\textbf{Case 1.} If $m\in(0,\frac{f(a+\eta(b,a))}{f(a)})$, then
	\begin{equation*}
	(s)\int_{a}^{a+\eta(b,a)} f(x) d\mu \leq min\{\beta,\eta(b,a)\},
	\end{equation*}
	where $\beta$ satisfies the following equation
	\begin{equation}\label{or38}
	(\eta(b,a)-\beta)^\alpha mf(\frac{a+\eta(b,a)}{m})-(\eta(b,a)-\beta)^\alpha f(a)-\eta(b,a)^\alpha (\beta-f(a))=0.
	\end{equation}
	\textbf{Case 2.} If $m=\frac{f(a+\eta(b,a))}{f(a)}$, then
	\begin{equation*}
	(s)\int_{a}^{a+\eta(b,a)} f(x) d\mu \leq min\{\beta,\eta(b,a)\},
	\end{equation*}
	where $\beta$ satisfies the following equation
	\begin{align}\label{or40}
	\nonumber & (\eta(b,a)-\beta)^\alpha \frac{f(a+\eta(b,a))}{f(a)}f(\frac{a+\eta(b,a)f(a)}{f(a+\eta(b,a))})\\
	&\quad-(\eta(b,a)-\beta)^\alpha f(a)-\eta(b,a)^\alpha (\beta-f(a))=0.
	\end{align}
	\textbf{Case 3.} If  $m\in (\frac{f(a+\eta(b,a))}{f(a)}, \eta(1,0))$, then
	\begin{equation*}
	(s)\int_{a}^{a+\eta(b,a)} f(x) d\mu \leq min\{\beta,\eta(b,a)\},
	\end{equation*}
	where $\beta$ satisfies the following equation
	\begin{align}\label{or42}
	\beta^\alpha(mf(\frac{a+\eta(b,a)}{m}))-\beta^\alpha f(a)- \eta(b,a)^\alpha(\beta-f(a))=0.
	\end{align}
\end{theorem}

Similarly we can prove the above results as Theorem \ref{of1}, \ref{of22}.

\end{large}
\end{document}